\documentclass[hidelinks]{article}

\usepackage{PRIMEarxiv}

\usepackage[utf8]{inputenc} 
\usepackage[T1]{fontenc}    
\usepackage{hyperref}       
\usepackage{url}            
\usepackage{booktabs}       
\usepackage{amsfonts}       
\usepackage{nicefrac}       
\usepackage{microtype}      
\usepackage{lipsum}
\usepackage{fancyhdr}       
\usepackage{graphicx}       
\graphicspath{{media/}}     

\usepackage{xspace}
\usepackage{fancybox}\usepackage{hyperref}
\usepackage{amsmath}
\usepackage{amsthm}
\usepackage{amssymb}
\usepackage{enumitem}

\title{Six equations in search of a finite-fold-ness proof
}

\author{
  Domenico Cantone \\
  Dept.\ of Mathematics and Computer Science \\
  University of Catania \\
  I-95125 Catania, Italy\\
  \texttt{domenico.cantone@unict.it} \\
   \And
  Luca  Cuzziol\\
  I-3100 Treviso, Italy\\
  \texttt{lucacuzz95@gmail.com} \\
  \AND
    Eugenio G. Omodeo\thanks{Corresponding author.} \\
  Dept.\ of Mathematics and Earth Sciences\\
  University of Trieste \\
  I-34127 Trieste, Italy\\
  \texttt{eomodeo@units.it} \\
}

\newcommand{\scnats}{\mbox{\scriptsize$\mathbb{N}$}}
\newcommand{\nats}{\mbox{$\mathbb{N}$}\xspace}
\newcommand{\ints}{\mbox{$\mathbb{Z}$}\xspace}
\newcommand{\rats}{\mathbb{Q}}
\newcommand{\tup}[1]{\left\langle{#1}\right\rangle}
\newcommand{\myvee}{\mbox{\it\small$\vee$}}
\newcommand{\mywedge}{\mbox{\it\small\&}}
\newcommand{\mybigwedge}{\mbox{\it\large\&}\,}
\newcommand{\myimp}{\:\mbox{\it\footnotesize$\Longrightarrow$}\:}
\newcommand{\mybimp}{\:\mbox{\it\footnotesize$\Longleftrightarrow$}\:}
\newcommand{\defFoAs}{\:\stackrel{\mbox{\tiny{Def}}}{\Longleftrightarrow}\:}
\newcommand{\defAs}{\:\stackrel{\mbox{\tiny{Def}}}{=}\:}
\newcommand{\M}[1]{\mbox{\boldmath$#1$}}
\newcommand{\mybox}[1]{\doublebox{\begin{minipage}{11.0cm}{\small #1}\end{minipage}}}
\newcommand{\smallmybox}[1]{\doublebox{\begin{minipage}{10.5cm}{\small #1}\end{minipage}}}
\newcommand{\overbar}[1]{\mkern 1.5mu\overline{\mkern-1.5mu#1\mkern-1.5mu}\mkern 1.5mu}
\newcommand{\pag}[1]{p.$\:${#1}}

\newcommand{\Sec}[1]{Sect.$\:$\ref{#1}}
\newcommand{\Secs}[2]{Sections$\:$\ref{#1} and \ref{#2}}
\newcommand{\abs}[1]{\vert #1 \vert}
\newcommand{\hug}[1]{\left\{#1\right\}}
\newcommand{\KOMMENT}[1]{}
\renewcommand{\qed}{\hfill{\boldmath$\dashv$}}
\newcommand{\Lem}[1]{Lemma~\ref{#1}}
\newcommand{\Thm}[1]{Theorem~\ref{#1}}
\newcommand{\Cor}[1]{Cor.~\ref{#1}}
\newcommand{\Fact}[1]{Fact~\ref{#1}}
\newcommand{\wt}{\,\mbox{\tt:}\:}
\newcommand{\Exa}[1]{Example~\ref{#1}}

\newtheorem{theorem}{Theorem}%
\newtheorem{lemma}{Lemma}%
\newtheorem{corollary}{Corollary}%
\newtheorem{example}{Example}%
\newtheorem{fact}{Fact}
\newtheorem{remark}{Remark}%

\begin{document}
\maketitle

\begin{abstract}
By following the same construction pattern which Martin Davis proposed in a 1968 paper of his, we have obtained six quaternary quartic Diophantine equations that candidate as `rule-them-all' equations: proving that one of them has only a finite number of integer solutions would suffice to ensure that each recursively enumerable set admits a finite-fold polynomial Diophantine representation.
\end{abstract}


\textbf{MSC Classification 2010/2020:} 03D25 , 11D25

\section*{Introduction}
\begin{sloppypar}\noindent In 1968, Martin Davis argued that Hilbert's 10$^{\rm th}$ problem would turn out to be algorithmically unsolvable if the equation
$\begin{array}{rcl}9\cdot\left(u^2+7\,v^2\right)^2-7\cdot\left(r^2+7\,s^2\right)^2&=&2\end{array}$
could be shown to admit only one solution in \nats (see \cite{Davis68}). Today we know that, more generally  (cf. \cite{DBLP:conf/cilc/CantoneO18,CCFO21,CO21}), if any of the equations\end{sloppypar}
\begin{center}\begin{tabular}{r||c}
Number $d$ & Associated quaternary quartic equation \\\hline
$\phantom{\overline{\overline{\mid}}}$2 & $\begin{array}{rcl}2\cdot\left({r}^{2} + 2\: s^{2}\right)^{2} - \left({u}^{2} + 2\: v^{2}\right)^{2}
&=&1\end{array}$\\
3 & $\begin{array}{rcl}3\cdot \left(r^2 +3\,s^2\right)^2 -\left(u^2 +3\,v^2\right)^2 &=& 2\end{array}$\\
7 & $\begin{array}{rcl}7\cdot\left(r^2+7\,s^2\right)^2-3^2\cdot\left(u^2+7\,v^2\right)^2&=&-2\end{array}$\\
11 & $\begin{array}{rcl}11\cdot\left( r^2+ r\: s\,+3\, s^2\right)^2-\left( v^2+ v\, u+3\, u^2\right)^2 &=& 2\end{array}$\\
19 & $\begin{array}{rcl}
19\cdot 3^2\cdot\left( r^2+ r\: s\,+5\, s^2\right)^2-13^2\cdot\left( v^2+ v\, u+5\, u^2\right)^2 &=& 2
\end{array}$\\
43 & $\begin{array}{rcl}43\cdot\left( r^2+ r\: s\,+11\, s^2\right)^2-\left( v^2+ v\, u+11\, u^2\right)^2 &=& 2
	\end{array}$, 
\end{tabular}\end{center}
associated with the respective Pell equations $x^2-d\:y^2=1$
 turned out to admit only a finite number of solutions in \ints, then every
recursively enumerable (`r.e.') subset of $\nats^{m}$---first and foremost the set of all triples $\tup{b,n,c}$ such that $b^n=c$---would admit a finite-fold polynomial Diophantine representation.

This means that to any partially computable function $f(\vec{a})$ from $\nats^m$ to \nats, there would correspond a multi-variate polynomial $D\in\ints[a_1,\dots,a_m,x_1,\dots,x_\kappa]$ such that for each
$\vec{a}\in\nats^m$ the following two conditions (representation and finite-fold-ness) hold:
\[\begin{array}{rclcl}i)&\phantom{xxxxxx}&f(\vec{a}) \mbox{ yields a value}&\mybimp&
\exists\,\vec{x}\in\nats^\kappa\ D(\vec{a},\vec{x})=0\,;\\
ii)&&\multicolumn{3}{l}{\exists\,b\in\nats\ \forall\,\vec{x}\in\nats^\kappa\;\big(\,D(\vec{a},\vec{x})=0\ \ \myimp\ \ b>\sum \vec{x}\,\big)\;.}
\end{array}\]
This would be the case if, say, the first of the above-listed quartic equations admitted only
the trivial solutions $r=\pm1$, $s=0$, $u=\pm1$, $v=0$;\footnote{In general, we call \emph{trivial} the solutions of the quartic associated with $d$ that correspond to the trivial solution $\tup{x,y}=\tup{1,0}$ of the Pell equation $x^2-d\,y^2=1$ (cf.\ the proof of \Lem{cor:three} for $d=19$). For every $d\neq43$ in the above table, these are the solutions with $r=\pm1$ and $s=0$. The case $d=43$ is exceptional: at $\tup{x,y}=\tup{1,0}$ the two factors $9\,x+59\,y$ and $59\,x+387\,y$ of $y_{2\ell+1}$ (see \Sec{sec:promising}) take the values $9$ and $59$, so that the trivial solutions are those with $r^2+r\,s+11\,s^2=9$ and $v^2+v\,u+11\,u^2=59$, e.g., $\tup{r,s,v,u}=\tup{3,0,3,2}$. Accordingly, for $d=43$ the non-triviality condition reads $r^2+r\,s+11\,s^2\neq9$, not $r\neq\pm1\ \myvee\ s\neq0$.} unfortunately, as will be reported in \Sec{sec:postscript}, this is not the case. 

\medskip

Why would it be important to establish whether any of the above equations has
only a finite number of solutions? The whole point is that if the equation associated with $d$ is finite-fold, then the following dyadic relation $\mathcal{J}_d$ over $\nats$ admits a polynomial Diophantine representation:
\[\begin{array}{rcrcl}
d\in\{2,7\}:&& \mathcal{J}_d(p,q) &\defFoAs& \exists\,\ell>4\,\Big[\,
q=\M{y}_{2^{\ell}}(d)\ \mywedge\ p\mid q\ \mywedge\ p\geq 2^{\ell+1}
\,\Big]\,,\\[0.14cm]
d\in\{3,11,19,43\}:
 && \mathcal{J}_{d}(p,q) &\defFoAs& \exists\,\ell>5\,\Big[\,
q=\M{y}_{2^{2\,\ell+1}}(d)\ \mywedge\ p\mid q\ \mywedge\ p\geq 2^{2\,\ell+2}\,\Big]\,,
\end{array}\]
where $\tup{\M{y}_i(d)}_{i\in\scnats}$ 
is the endless, strictly ascending, sequence consisting of all non-negative integer solutions to the above-mentioned Pell equation $d\,y^2+1=\Box$\,.\footnote{`$Q=\Box$' means that $Q$ must be a perfect square.} 
Independently of representability, each $\mathcal{J}_d$ turns out to satisfy Julia Robinson's exponential growth criteria, namely:

\begin{itemize}
\item $\mathcal{J}_d(p,q)$ implies $p>1\ \mywedge \ q\leq p^p\:$,

\vspace{.12cm}

\item for each $k\geq 0$, there are $p$ and $q$ such that $\mathcal{J}_d(p,q)$ and $p^k<q\:$;
\end{itemize}

\noindent hence, as by \cite{Robinson69a}, representability would propagate from $\mathcal{J}_d$ to any r.e.\ set.\footnote{Notice that removal of a finite number of pairs from a relation $\mathcal{J}\subseteq\nats\times\nats$ does not disrupt
its exponential-growth property nor its Diophantine representability (if $\mathcal{J}$ enjoys one or the other).}

\medskip

In \cite{Matiyasevich2010}, Yuri Matiyasevich noted that the property
\[\tag{\boldmath$\ddagger$}\label{propertyByMatiyasevich}
\smallmybox{integers $\alpha>1\,,\:\beta\geq0\,,\:\gamma\geq0\,,\delta>0$ exist such that to each $w\in\nats$, $w\neq 0$, there correspond $p,q$ such that $\mathcal{J}(p,q)$, $p<\gamma w^\beta$, and $q>\delta\alpha^w$ hold}
\] is such that if an exponential-growth relation $\mathcal{J}\subseteq\nats\times\nats$ enjoying it admits a polynomial Diophantine \emph{finite-fold} representation, then exponentiation---and, consequently, any r.e.\ set---admits, in turn, a representation
of that kind. 
We have verified that each one of our $\mathcal{J}_d$'s actually enjoys this property, by fixing $\alpha,\beta,\gamma$, and $\delta$
as follows:

\[\begin{array}{rcccccccl}
d\in\{2,7\}:&& \alpha=\M{x}_1(d)\,,&&\beta=2\,,&&\gamma=2^7\,,&&\delta=1\,,\\[0.14cm]
d\in\{3,11,19,43\}:
 && \alpha=\big(\M{x}_1(d)\big)^4\,,&&\beta=2\,,&&\gamma=2^{15}\,,&&\delta=\big(\M{x}_1(d)\big)^3\,,
\end{array}\]

\noindent where $\M{x}_1(d)$ is the positive integer such that $d\cdot\big(\M{y}_1(d)\big)^2+1=\big(\M{x}_1(d)\big)^2$\,.

\medskip

At this point, in order to complete our `quest for Diophantine finite-fold-ness', we should prove that at least one of our quaternary quartics has, in all, finitely many solutions in \ints. Contrary to Davis's original expectation, Daniel Shanks
and Samuel S. Wagstaff conjecture that his quaternary quartic
(the one associated with $d=7$) has infinitely many solutions; in \cite{SW95}, they proposed a method that helps in finding  non-trivial solutions to equations of that kind, and they exploited it to actually detect some 50 solutions to Davis's equation.

The Shanks--Wagstaff method also brought to light non-trivial solutions to the equations associated with $d=3$\,, $d=11$\,, and $d=2$\,, but
it turned out impractical to use it in connection with $d=19$ and $d=43\,$, whose treatment calls for the factorization of  huge numbers into primes. (Added in revision: a non-trivial solution to the equation associated with $d=43$ was nevertheless found by L. Cuzziol and P. Campochiaro, namely $r=2644220400$, $s=-305590179873$, $v=2488975045142$, $u=134875386175$; it corresponds to the Pell solution $\tup{x_6,y_6}$ of $x^2-43\,y^2=1$.)

\begin{center}------------\end{center}
As said at the beginning, this paper prolongs a triad of papers putting in context, revisiting, and extending Davis's construction of a potential rule-them-all equation
(see  \cite{Davis68}). The papers \cite{Mat1974,DMR76}
saw in Davis's approach a strategy to affirmatively answer the question: ``Is there a finitefold (or better a singlefold) polynomial Diophantine definition of $\,b^n=c$\,?''.  Recently, the interest in that finite-foldness issue was aroused again in \cite{Mat06,Matiyasevich2010}, and \cite{Mat016} reiterated the relevance of Davis's approach for solving it. 

A recapitulation of the stimuli for revisiting \cite{Davis68} today can be found in \cite{CCFO21}, which also offers a bird's-eye view of how to treat the number $d=2$. The number $d=11$ has been treated in \cite{DBLP:conf/cilc/CantoneO18}; the number $d=3$ has been treated in great detail in \cite{CO21}, and this exempts us---we deem---from providing many details concerning the treatment of $d=19$, to which we devote a biggish part of \Sec{sec:promising}.

In \Sec{sec:propertyByMatiyasevich} we produce a proof that property \eqref{propertyByMatiyasevich} holds for the paradigmatic case $\mathcal{J}=\mathcal{J}_{19}$\,; we also indicate how to shift this proof to the case $\mathcal{J}=\mathcal{J}_{2}$\,. 
In \Secs{sec:attemptToFindNonTrivialFor2}{sec:postscript} we illustrate how to look systematically for a non-trivial solution to one of the quaternary quartic equations under study and report on successful findings based on this approach.

\Sec{sec:uniqueFactorizationRings} surveys some key number-theoretic prerequisites.

\section{Unique-factorization rings of the integers of $\rats(\sqrt{-d}\,)$}\label{sec:uniqueFactorizationRings}

Davis's construction of a potential rule-them-all equation exploits a square-free rational integer $d>0$ such that in the imaginary quadratic field $\rats(\sqrt{-d}\,)$ the algebraic integers\footnote{The algebraic integers mentioned here are the elements of $\rats(\sqrt{-d}\,)$ which are roots of monic polynomials $x^n + c_{n-1}\,x^{n-1} + \cdots + c_0$ with rational integer coefficients $c_i\,$.}  form a unique-factorization integral domain $\mathcal{A}_d$\,. All such  numbers were discovered by Carl Friedrich Gauss; they are \\\centerline{1, 2, 3, 7, 11, 19, 43, 67, 163 (see \href{https://oeis.org/A003173}{OEISA003173}),} and they are often termed {\bf\em Heegner numbers} after the name of the scholar who gave, in the 1950s, a decisive contribution to the proof that no more numbers with the desired property exist (cf. \cite{Hee52,Sta69}). Number $1$ must be discarded beforehand, because it is a perfect square (hence $x^2-1\,y^2$ is not a Pell's equation); as it turns out that $d\ \equiv\ 3\ \mod 4$
holds for all other Heegner numbers except $d=2$, we have: 
        \[\begin{array}{rcl}\mathcal{A}_d&=&\left\{\begin{array}{rl}
        \ints\left[\,\sqrt{-d}\,\right]&\mbox{when $d=2$},\\[0.25cm]
        \ints\left[\,\frac{1+\sqrt{-d}}{2}\,\right]&\mbox{for $d\in\{3, 7, 11, 19, 43, 67, 163\}$}.\end{array}\right.\end{array}\]

The rational prime numbers that cease to be prime in this ring $\mathcal{A}_d$ turn out to be the ones writable in the following {\bf\em norm form} (where $w,t\in\ints$):

\[\begin{array}{lcl}
w^2+d\,t^2 && \mbox{when $d=2$,}\\
w^2+w\,t+\frac{d+1}{4}\,t^2  &\phantom{xx}& \mbox{for $d\in\{3, 7, 11, 19, 43, 67, 163\}$.}
\end{array}\]

\noindent(Remark: In the case when $d=3$, the numbers of this form coincide with the ones writable in the form $w^2+3\,t^2$; 
when $d=7$, then $2$ is the sole prime number of norm form that cannot be written in the form $w^2+7\,t^2$).

For any Heegner number $d\neq 1$, we call {\bf\em inert} those numbers $p$, prime in \ints, which remain prime in the enlarged ring $\mathcal{A}_d$\,. Notice that it is $\mathcal{A}_d$, not $\ints\left[\sqrt{-d}\,\right]$, that matters here; e.g., 
$p=2$ is irreducible, but not prime, in $\ints\left[\sqrt{-3}\,\right]$\,;\footnote{To see that 2 is not prime in $\ints\left[\sqrt{-3}\,\right]$\,, it suffices to observe that 2 divides $4=(1+\sqrt{-3})\,(1-\sqrt{-3})$ without dividing either factor; on the other hand $2$ is irreducible in $\ints\left[\sqrt{-3}\,\right]$\,, in the sense it can be factorized as $2=(a+b\,\sqrt{-3})\,(a'+b'\,\sqrt{-3})$ only with $b=b'=0$ and $a=\pm1\vee a'=\pm1$.} however, it becomes prime in the ring of integers of $\rats(\sqrt{-3}\,)$, which is larger than $\ints\left[\sqrt{-3}\,\right]$\,. \hfill The number $d$ itself is \emph{not} an inert prime.

\noindent Relative to a $d\neq 2$\,, an odd prime $p$ other than $d$ turns out to be \emph{inert} if and only if 

\[\begin{array}{lcccll}-d&\not\equiv& x^2 \mod p &&& \mbox{holds for every $x\in \ints$ such that $\abs{x}<p$\,;}\end{array}\]

\noindent {thanks to the quadratic reciprocity law (see \cite{cox2013primes}), this also amounts to the property

\[\begin{array}{lcccll}p&\not\equiv& y^2 \mod d &&& \mbox{for every $y\in \ints$ such that $0<\abs{y}<d$\,,}\end{array}\]}
which offers a practical criterion to recognize inert primes and remains valid also for $p=2$.

On these grounds one easily preps the table below, showing, relative to each Heegner number $d\notin\hug{1,67,163}$, which primes are inert:\footnote{Including $d=67$ and $d=163$ in this table of inert primes would have been
easy but pointless, since we must refrain from fully treating those two numbers in this paper.}
\begin{center}\footnotesize
\begin{tabular}{r||cl|cl}
$d$ && \hspace{0.2cm}{inert prime $p$} && \hspace{0.2cm}{representable prime $q$ other than} $d$ \\ \hline
2 && $p \equiv 5,7 \mod 8$ && $q \equiv 1,3 \mod 8$\\
3 && $p \equiv 2 \mod 3$ && $q \equiv 1 \mod 3$\\
7 && $p \equiv 3,5,6 \mod 7$ && $q \equiv 1,2,4 \mod 7$ \\
11 && $p \equiv 2,6,7,8,10 \mod 11$ && $q \equiv 1,3,4,5,9 \mod 11$ \\
19 && $p \equiv 2,3,8,10,12,13,14,15,18 \mod 19$ && $q \equiv 1,4,5,6,7,9,11,16,17 \mod 19$\\
43 && $p \equiv 2,3,5,7,8,12,18,19 \mod 43$ && $q \equiv 1,4,6,9,10,11,13,14 \mod 43$\\
''\:&& $p \equiv 20,22,26,27,28,29 \mod 43$ && $q \equiv 15,16,17,21,23,24,25 \mod 43$\\
''\:&& $p \equiv 30,32,33,34,37,39,42 \mod 43$ && $q \equiv 31,35,36,38,40,41 \mod 43$\KOMMENT{\\
67 & & \small$\equiv 1,4,6,9,10,14,15,17 \mod 67$ \\
$''$ & & \small$\equiv 19,21,22,23,24,25 \mod 67$ \\
$''$ & & \small$\equiv 26,29,33,35,36,37 \mod 67$ \\
$''$ & & \small$\equiv 39,40,47,49,54,55 \mod 67$ \\
$''$ & & \small$\equiv 56,59,60,62,64,65 \mod 67$ \\
163 & & \small$\equiv 1,4,6,9,10,14,15,16 \mod 163$ \\
$''$ & & \small$\equiv 21,22,24,25,26,33 \mod 163$ \\
$''$ & & \small$\equiv 34,35,36,38,39,40 \mod 163$ \\
$''$ & & \small$\equiv 41,43,46,47,49,51 \mod 163$ \\
$''$ & & \small$\equiv 53,54,55,56,57,58 \mod 163$ \\
$''$ & & \small$\equiv 60,61,62,64,65,69 \mod 163$ \\
$''$ & & \small$\equiv 71,74,77,81,83,84 \mod 163$ \\
$''$ & & \small$\equiv 85,87,88,90,91,93 \mod 163$ \\
$''$ & & \small$\equiv 95,96,97,100,104 \mod 163$ \\
$''$ & & \small$\equiv 111,113,115,118,119 \mod 163$ \\
$''$ & & \small$\equiv 121,126,131,132,133 \mod 163$ \\
$''$ & & \small$\equiv 134,135,136,140,143 \mod 163$ \\
$''$ & & \small$\equiv 144, 145,146,150,151 \mod 163$ \\
$''$ & & \small$\equiv 152, 155,156,158,160 \mod 163$ \\
$''$ & & \small$\equiv 161 \mod 163$}
\end{tabular}
\end{center}

In \cite{Davis68}, inert primes are called \emph{poison primes}, because they `poison' those $m\in\nats$ in whose standard  factorization they appear with an odd exponent, in the following sense:
if $m$ gets so poisoned, then it is not representable in the  above-specified quadratic form. We sum up the situation in the frame shown here below:

\vspace{0.15cm}

\begin{center}
\mybox{{\bf\em Inert}, or `{\em poison}' prime:\begin{itemize}\item it remains prime;
\item it is not representable in the norm form;
\item it poisons (by rendering it, in turn, not so representable) any $m\in\nats$ in whose standard factorization it occurs with an odd exponent. \end{itemize}}
\end{center}

\begin{example}\label{exa:reprMod2}
Take $d=2$. From the discussion just given, it follows that any odd representable number $m$ satisfies either one of the congruences $m\equiv d\pm1\mod 8$\,. \qed
\end{example}

\section{Promising (Diophantine?) exponential-growth relations} \label{sec:promising}
In this section we outline Davis's construction of a potential rule-them-all equation,
choosing the Heegner number $d=19$ as our running example. Before going into technicalities related to that number, we offer some background information.

\medskip

To each Heegner number $d\neq 1$ there corresponds the Pell equation
\[\begin{array}{rcl}d\,y^2+1=x^2\,,\end{array}\]

\noindent trivially solved by $y_0=0\,,\:x_0=1$\,. Its fundamental solution  $y_1\,,\:x_1$ (namely, the one with $y,x$ positive integers of smallest possible size), and the rule for getting from two consecutive solutions
$y_k\,,\:x_k$ and $y_{k+1}\,,\:x_{k+1}$ its next solution $y_{k+2}\,,\:x_{k+2}$\,, are tabulated here below:
\begin{center}\small\begin{tabular}{r||r|r||c|c}
$d$ & $y_1$ & $x_1$ & $y_{k+2}$ & $x_{k+2}$ \\\hline
2 & 2 & 3 & $\;\;\;\ \ \ \ 6\,y_{k+1}-y_k$ & $\;\;\;\ \ \ \ 6\,x_{k+1}-x_k$\\
3 & 1 & 2 & $\;\;\;\ \ \ \ 4\,y_{k+1}-y_k$ & $\;\;\;\ \ \ \ 4\,x_{k+1}-x_k$\\
7 & 3 & 8 & $\;\;\:\ \ \ 16\,y_{k+1}-y_k$  & $\;\;\:\ \ \ 16\,x_{k+1}-x_k$\\
11 & 3 & 10 & $\;\;\,\ \ \ 20\,y_{k+1}-y_k$ & $\;\;\,\ \ \ 20\,x_{k+1}-x_k$\\
19 & 39 & 170 & $\;\:\,\ \ 340\,y_{k+1}-y_k$ & $\;\:\,\ \ 340\,x_{k+1}-x_k$\\
43 & 531 & 3482 & $\;\,\ 6964\,y_{k+1}-y_k$  & $\;\,\ 6964\,x_{k+1}-x_k$\\
67 & 5967 & 48842 & $\; 97684\,y_{k+1}-y_k$ & $\; 97684\,x_{k+1}-x_k$\\
163 & 5019135 & 64080026 & $2\cdot x_1\cdot y_{k+1}-y_k$ & $2\cdot x_1\cdot x_{k+1}-x_k$
\end{tabular}\end{center}

\bigskip

A cornerstone in the construction of candidate rule-them-all equations is the following basic fact concerning Pell equations (cf., e.g., \cite[Corollary~3.2]{CO21}): 
\begin{lemma}\label{lem:x_twoell}
Consider the Pell equation $x^2-\delta\,y^2=1$, with $\delta>0$ a non-square integer; let
$\tup{x_0,y_0},\:\tup{x_1,y_1},\:\tup{x_2,y_2},\dots$ be its solutions in \nats, with $y_0 < y_1 < y_2<\cdots$\,. Then the following identities hold for every $\ell\in\nats$\,: \\
\[\begin{array}{rclcrcl}x_{2\,\ell}&=&x_\ell^2+\delta\,y_\ell^2&\ \ \mbox{\em and}\ \ &
y_{2\,\ell}&=&2\,x_\ell\,y_\ell\,.\end{array}\]
\end{lemma}

\vspace{-0.7cm}

$\phantom{}$\qed

\bigskip

\noindent This lemma shows that each $x_{2\,\ell}$ entering in a representation of $1$ in the form $x^2-\delta\,y^2$ is, in turn, represented by the quadratic form $x^2+\delta\,y^2$\,. One can get \Lem{lem:x_twoell} from the equality $x_k+y_k\sqrt{\delta}=\big(x_{1}+y_{1}\sqrt{\delta}\,\big)^k$ holding for all $k\in\nats$\,,
by taking into account the irrationality of $\sqrt{\delta}$\,. Consequently, since each entry of the sequence $\big\langle x_{2\,\ell+1} + y_{2\,\ell+1}\sqrt{\delta}\big\rangle_{\mbox{\scriptsize$\ell\in\nats$}}$ equals $\big(x_{2\,\ell}+y_{2\,\ell}\sqrt{\delta}\,\big)\:\big(x_{1}+y_{1}\sqrt{\delta}\,\big)$, we have
\[\begin{array}{rclcrcl}y_{2\,\ell+1}&=&y_1\,x_\ell^2+2\,x_1\,x_\ell\,y_\ell+\delta\,y_1\,y_\ell^2\:,\end{array}\]
and so, when $\delta$ is one of the $d$'s that interest us here:
\begin{center}\small\begin{tabular}{rc||cl}
$d=2$ &&& $y_{2\,\ell+1}=2\,(x_\ell+y_\ell)\:(x_\ell+2\,y_\ell)$ \\
$d=3$ &&& $y_{2\,\ell+1}=(x_\ell+y_\ell)\:(x_\ell+3\,y_\ell)$ \\
$d\in\hug{7\,,\:11}$&&& $y_{2\,\ell+1}=(x_\ell+3\,y_\ell)\:(3\,x_\ell+d\:y_\ell)$ \\
$d=19$ &&& $y_{2\,\ell+1}=(3\,x_\ell+13\,y_\ell)\:(13\,x_\ell+3\cdot 19\,y_\ell)$ \\
$d=43$ &&& $y_{2\,\ell+1}=(9\,x_\ell+59\,y_\ell)\:(59\,x_\ell+9\cdot 43\,y_\ell)$ \\
$d=67$ &&& $y_{2\,\ell+1}=(27\,x_\ell+221\,y_\ell)\:(221\,x_\ell+27\cdot67\,y_\ell)$ \\
$d=163$ &&& $y_{2\,\ell+1}=(627\,x_\ell+8005\,y_\ell)\:(8005\,x_\ell+627\cdot163\,y_\ell)$ \\
\end{tabular}\end{center}

It turns out that the binomials appearing in each one of the above equalities
evaluate to co-prime numbers  $v_\ell\,,\:w_\ell$; i.e.: $\gcd(x_\ell+y_\ell,x_\ell+2\,y_\ell)=1$ when $d=2$, $\gcd(x_\ell+y_\ell,x_\ell+3\,y_\ell)=1$ when $d=3$, \dots,
$\gcd(627\,x_\ell+8005\,y_\ell,8005\,x_\ell+102201\,y_\ell)=1$ when $d=163$\,. 
Accordingly, if $y_{2\,\ell+1}$ is representable in the norm form, so are its factors  $v_\ell$ and $w_\ell$; in fact, no inert prime could possibly divide either $v_\ell$ or $w_\ell$ to an odd power, else it would divide $y_{2\,\ell+1}$ to the same odd power.\footnote{The consideration made in this paragraph may need some
refinement for particular $d$'s: this is apparent in the treatment of $d=7$ in \cite{Davis68}, and will surface again in the treatment of $d=19$ inside the proof of \Thm{thm:repr19} below.}

\begin{example}
Take $d=19$. To see that $\gcd(3\,x_\ell+13\,y_\ell,13\,x_\ell+57\,y_\ell)=1$ holds for every $\ell$\,,
begin by observing that $\gcd(x_\ell,y_\ell)=1$\,: in fact,
any positive integer $t$ dividing both of $x_\ell$ and $y_\ell$
will divide $x_\ell^2-d\,y_\ell^2$, which equals $1$; therefore such a $t$
must equal $1$\,. Next, by way of contradiction, suppose that a prime number $p$ exists such that
$p\mid 13\,x_\ell+57\,y_\ell$ and $p\mid 3\,x_\ell+13\,y_\ell$\,. Then $p\mid 2\,y_\ell$ holds, because $2\,y_\ell=3\,(13\,x_\ell+57\,y_\ell)-13\,(3\,x_\ell+13\,y_\ell)$\,; therefore either $p=2$ or $p\mid y_\ell$\,. Inspection of the Pell equation
at hand, whose discriminant $d$ is odd, shows us that
$y_\ell$ and $x_\ell$ have opposite parity; so $ 3\,x_\ell+13\,y_\ell$ is odd and $p\neq 2$\,. Hence $p\mid y_\ell$, and moreover $p\mid [9\,(3\,x_\ell+13\,y_\ell)-2\,(13\,x_\ell+57\,y_\ell)]$\,, i.e. $p\mid x_\ell+3\,y_\ell$\,, and so $p\mid x_\ell$\,, which contradicts the co-primality between $x_\ell$ and $y_\ell$\,.\footnote{For each odd $d$, the analogous co-primality result is proved very much like in this example, with only the final step requiring a bit of ingenuity. When $d=3$, the clue is that $2\,(x_\ell+3\,y_\ell)-(x_\ell+y_\ell)=x_\ell-5\,y_\ell$. When $d=7$ or $d=11$, the clue is that $(x_\ell+3\,y_\ell)-3\,y_\ell=x_\ell$; when $d=43$, it is that $2\,(59\,x_\ell+387\,y_\ell)-13\,(9\,x_\ell+59\,y_\ell)=x_\ell+7\,y_\ell$; when $d=67$, the clue is that $11\,(221\,x_\ell+1809\,y_\ell)-90\,(27\,x_\ell+221\,y_\ell)=x_\ell+9\,y_\ell$; when $d=163$, it is that $73\,(8005\,x_\ell+102201\,y_\ell)-932\,(627\,x_\ell+8005\,y_\ell)=x_\ell+13\,y_\ell$.

When $d=2$, the contradiction $p\mid y_\ell\ \mywedge\ p\mid x_\ell$ follows from supposing that a prime number $p$ exists such that $p\mid x_\ell+y_\ell$ and $p\mid x_\ell+2\,y_\ell$: in fact, $y_\ell=(x_\ell+2\,y_\ell)-(x_\ell+y_\ell)$ and $x_\ell=2\,(x_\ell+y_\ell)-(x_\ell+2\,y_\ell)$\,.
}
\qed\end{example}

What follows will adjust to the case $d=19$ a general technique for associating with a given Heegner number $d\neq 1$ an equation whose solutions correspond to the representations of a certain integer $C$ by a quadratic form $d\cdot a^2\cdot X^2 - b^2\cdot Y^2$, where $X,Y$ are in turn representable by the quadratic norm form associated with $d$ as seen in \Sec{sec:uniqueFactorizationRings}. In most cases, either $a=1$ or $b=1$ or $a=b=1$ holds; $C$ often equals 2. Inspection of the nearly complete table of resulting quartic equations shown at the beginning of the Introduction makes it evident that the quaternary quartic equation associated with $d=19$ is relatively elaborate; thus, our running example should encompass all technical difficulties. 

\subsection{Construction of the quaternary quartic associated with $d=19$}

Two
lemmas will aid in the proof of a crucial statement to be seen below.
\begin{lemma} \label{lem:repr19a} Take $d=19$. For $m,h=1,2,3\dots$
it holds that \[\begin{array}{rcl}y_{2^m\cdot h}&=&2^m \,x_h \,y_h\cdot\textstyle\prod_{0<i<m}x_{2^i\cdot h}\,.\end{array}\]
In particular, we have $\begin{array}{rcl}y_{2^m}&=&2^{m+1} \cdot3\cdot5\cdot13\cdot17\cdot\textstyle\prod_{0<i<m}x_{2^i}\end{array}$\,. 
\end{lemma}
\begin{proof}
The claim is proved by induction on $m$\,: it readily follows from \Lem{lem:x_twoell} when $m=1$\,; moreover, when $m=k+2$\,, \Lem{lem:x_twoell} together with the induction hypothesis yields that
\[\begin{array}{rcccccl}y_{2^{m}\cdot h}&=&y_{2^{k+2}\cdot h}&=&2\,x_{2^{k+1}\cdot h}\,y_{2^{k+1}\cdot h}&=&
2^{k+2}\,x_h\,y_h\cdot \textstyle\prod_{0<i\leq k+1}x_{2^i\cdot h}\\[0.12cm]&&&&&=&2^{m}\,x_h\,y_h\cdot \textstyle\prod_{0<i<m}x_{2^i\cdot h}\,.\end{array}\]

\vspace{-1cm}

$\phantom{}$

\end{proof}

Note that the recursive rule for calculating the $y_k$'s yields,
since $y_1=39$, that $39\mid y_k$ for every $k$\,.
\begin{lemma} \label{lem:repr19}
Take $d=19$. If ${y_n}/{39}$ is representable, where $n=2^m\cdot h$, $h$ is odd, and $m>0$, then ${y_h}/{39}$ is representable.
\end{lemma}
\begin{proof} By way of contradiction, suppose that $n=2^m\cdot h$, $h$ is odd, $m>0$, but an inert prime $p$ exists dividing $y_h/39$ to an odd power. Since $h$ is odd, $y_h$ is odd and $p\neq 2$. We know from \Lem{lem:x_twoell} that $x_{2^i\cdot h}$ is representable for each $i>0$: in fact $x_{2^i\cdot h}$ can be written in the form $u^2+19\,v^2$\,, hence it can also be written as $(u-v)^2+(u-v)(2\,v)+5\,(2\,v)^2$, and hence it has the norm form pertaining to $d=19$\,.	Therefore $p$ divides $x_{2^i\cdot h}$ to
	an even power (perhaps $0$), for $i=1,\dots,m-1$\,.
	In addition we have $p\ \nmid \ 2^{m}$; moreover, $p\ \nmid \ x_h$, because $\gcd(x_h,y_h)=1$.
	So we finally use \Lem{lem:repr19a} to prove that $p$ divides $y_n/39$ to an odd power, thus 
	getting the sought contradiction.
\end{proof}

Here is the above-mentioned key statement:
\begin{theorem} \label{thm:repr19}
Take $d=19$. If $y_n/39$ is representable (in the norm form $w^2+w\,t+5\,t^2$) for some $n>0$ not a power of $2$, then the constraint system
\[\left\{\begin{array}{lcrlcl}
X^2&-&19\cdot 39^2\,\cdot& Y^2 &=& 1\,,\\
X&+&13^2\,\cdot& Y&=&r^2+r\,s+5\:s^2\,,\\
X&+&19\cdot 3^2\,\cdot& Y&=&v^2+v\,u+5\:u^2\,,\\
&&&Y&>&0
\end{array}\right.\]
has a solution $\bar X,\bar Y,\bar r,\bar s,\bar v,\bar u$ in \ints such that either $\bar r\neq \pm 1$ or $\bar s\neq 0$ holds
and, moreover, $39\cdot(\bar X+13^2 \cdot \bar Y)\cdot(\bar X+19\cdot 3^2\cdot\bar Y)\mid y_n$\,.
\end{theorem}
\proof{Suppose $n=2^m(2\,\ell+1)$, with $m\geq0$ and $\ell>0$ (the latter because $n$ is not a power of $2$), be such that $y_n/39$ is representable. Put $h=2\,\ell+1$. If $m>0$, \Lem{lem:repr19} tells us that $y_h/39$ is representable; if $m=0$, this is our very hypothesis. Thus $y_{2\,\ell+1}/39$ is representable; therefore,
$x_\ell+3\cdot19\frac{y_\ell}{13}$ and $x_\ell+13\,\frac{y_\ell}{3}$ are also representable. By setting $\bar X=x_\ell$ and 
$\bar Y=\frac{y_\ell}{39}$, we hence have numbers $\bar r,\bar s,\bar v,\bar u$ that satisfy the system appearing in the claim. 

It is untenable that $\bar r =\pm1\ \mywedge\ \bar s=0$\,, because this would imply $x_\ell+13^2\,\frac{y_\ell}{39}=1$\,, contradicting $x_\ell>y_\ell\geq 39$\,.

\Lem{lem:repr19a} tells us that $y_h\mid y_{2^m\cdot h}$ (trivially so when $m=0$),
whence, in the present context:
$39\:(\bar X+13^2 \: \bar Y)\:(\bar X+19\cdot 3^2\:\bar Y)=(13\,x_\ell+3\cdot 19\,y_\ell)\:(3\,x_\ell+13\,y_\ell)=y_{2\,\ell+1}\mid y_n$\,. \qed}

\begin{corollary}\label{cor:repr19} Under the hypothesis of \Thm{thm:repr19}, the equation
\[\tag{\boldmath$\dagger$}\label{eq:rulMall19}\begin{array}{rcl}
19\cdot 3^2\cdot\left( r^2+ r\: s\,+5\, s^2\right)^2-13^2\cdot\left( v^2+ v\, u+5\, u^2\right)^2 &=& 2
\end{array}\]
has a solution $\bar r\,,\:\bar s\,,\:\bar v\,,\:\bar u$
such that either $\bar r\neq\pm1$ or $\bar s\neq0$ holds and, moreover, $39\:(\bar r^2+ \bar r\: \bar s\,+5\, \bar s^2)\:(\bar v^2+ \bar v\, \bar u+5\, \bar u^2)\mid y_n$\,.
\end{corollary}
\proof{
	Indeed, by solving the system of constraints of \Thm{thm:repr19} in the manner discussed there, we will have
	
\smallskip

	\centerline{\small$\begin{array}{cl}
	&19\cdot3^2\,\cdot\,\left(\,\bar r^2+ \bar r\, \bar s+5\, \bar s^2\,\right)^2-
	13^2\,\cdot\,\left(\,\bar v^2+ \bar v\, \bar u+5\, \bar u^2\,\right)^2\\
	=&171\cdot\left(\bar X^2+338\cdot \bar X\cdot \bar Y+169^2\cdot \bar Y^2\right)-169\cdot\left(\bar X^2+342\cdot \bar X\cdot \bar Y+171^2\cdot \bar Y^2\right)\\
	=&2\cdot\left(\bar X^2-19\cdot 39^2\cdot \bar Y^2\right)\\
	=&2\;,
	\end{array}$}

\smallskip

\noindent where the {\bf\em non-triviality condition} $\bar r\neq\pm1\ \myvee\ \bar s\neq 0$\,, along with the divisibility constraint stated in the claim, are satisfied. \qed
}

\subsection{Is $\left\{\,y_{2^{2\,\ell+1}}(19)\,:\:\ell=0,1,2,\dots\,\right\}$ a Diophantine set?}\label{sec:assertionH19}

We continue to refer to the Heegner number 19. Let $\mathcal{H}$ stand for the assertion: \begin{center}\mybox{\em\normalsize The equation \eqref{eq:rulMall19} admits, altogether, finitely many solutions in integers.}\end{center}

Then, by combining \Cor{cor:repr19} with the preparatory statement that $y_{2^m}/39$ is representable if and only if $m$ is odd
(a fact easily obtainable from the remark made at the end of the claim
of \Lem{lem:repr19a}), we get:
\begin{lemma}\label{cor:three} 
$\mathcal{H}$ implies that $\{\,y_{2^{2\ell+1}}\wt \ell\geq 0\,\}$ admits a polynomial Diophantine representation.
\end{lemma}
\begin{proof}
We know from \Cor{cor:repr19} that if $y_{n}/39$ is representable for some $n > 1$ not a power of 2, hence of the form $n = 2^{m}\,(2\,\ell+1)$ with $m \geq 0$ and $\ell > 0$, then the equation \eqref{eq:rulMall19} has a non-trivial integer solution
$\bar{r}, \bar{s}, \bar{v}, \bar{u}$ such that
$39\:(\bar r^2+ \bar r\: \bar s\,+5\, \bar s^2)\:(\bar v^2+ \bar v\, \bar u+5\, \bar u^2)\mid y_n$\,.

This, along with the preparatory lemma, yields the following \emph{sufficient} conditions in order for the property
\begin{equation}\label{DiophantinePred}
y \in \hug{y_{2^{2\,\ell+1}}\wt \ell \geq 0 }
\end{equation}
to hold:
\begin{enumerate}[label=(\roman*)]

\item\label{second} $y = y_{n}$, for some $n > 1$;
\smallskip
\item\label{third} $\frac{y}{39}$ is representable in the form $w^2\pm wt+5\,t^2$\, with $w,t\in\nats$;
\smallskip
\item\label{fourth} $39\:(\bar r^2+ \bar r\: \bar s\,+5\, \bar s^2)\:(\bar v^2+ \bar v\, \bar u+5\, \bar u^2) \nmid y$, for any solution $\bar r,\bar s,\bar v, \bar u$ to \eqref{eq:rulMall19} such that $\bar r\neq\pm1\ \myvee\ \bar s\neq 0$\,.
\end{enumerate}

Here are existential Diophantine definitions,  in \nats, of the first two of these:
\begin{description}
\item{\;\ref{second}} $\quad 19\,y^2+1=\Box\ \mywedge\ y>39$ (since $39=y_1$)\,;
\item{\ref{third}} $\exists\,w\:\exists\,t\:\left[\left(39\,w^2+5\cdot39\:t^2-y\right)^2 = 39^2\:w^2\,t^2 \right]$\,.
\end{description}

\noindent Moreover, if \eqref{eq:rulMall19} admits only finitely many solutions, then also \ref{fourth} admits an existential Diophantine definition. Indeed, let $\tup{r_{1},s_{1},v_{1},u_{1}},\dots,$ $\tup{r_{\kappa},s_{\kappa},v_{\kappa},u_{\kappa}}$ be all of the non-trivial solutions to \eqref{eq:rulMall19} that have $s_i,u_i\in\nats$ and $r_i,v_i\in\ints$ for each $i$. Then \ref{fourth} is easily seen to be statable over \nats as
\begin{multline*}
\!\!\!\!\!\!(\exists\, q_{1},\ldots,q_{\kappa},\,g_{1},\ldots,g_{\kappa},\, z_{1},\ldots,z_{\kappa})\\
 \mybigwedge_{i=1}^{^{\scriptstyle\kappa}}  \Big[\;y = 39\,\left( r_i^2+ r_i\: s_i\,+5\: s_i^2\right)\,\left( v_i^2+ v_i\, u_i+5\: u_i^2\right)\,q_i+g_i+1\ \mywedge \\
 {}g_{i}+z_{i}+2 = 39\,\left( r_i^2+ r_i\: s_i\,+5\: s_i^2\right)\,\left( v_i^2+ v_i\, u_i+5\: u_i^2\right)\Big]\,.
\end{multline*} 
 
In order to complete the proof that the membership relation \eqref{DiophantinePred} is Diophantine when \eqref{eq:rulMall19} admits only finitely many solutions, it only remains to be shown that \ref{second}--\ref{fourth} also are necessary conditions for \eqref{DiophantinePred} to hold. This will result in a polynomial Diophantine representation of the property
$y\in\{\,y_{2^{2\ell+1}}\,:\:\ell\geq 0\,\}$, \textbf{if} {\em the number of solutions to  the equation \eqref{eq:rulMall19} is finite}\,! (An issue that we are unable to answer.) 

\smallskip

Let, hence, $y = y_{2^{2\,\ell+1}}$ hold for some $\ell \geq 0$. Obviously \ref{second} holds and, by the preparatory lemma stated (without proof) in front of this one, \ref{third} holds as well.

Towards a contradiction, suppose that we have
\begin{equation}\label{contrHyp}
39\:\left(\bar r^2+ \bar r\: \bar s\,+5\, \bar s^2\right)\:\left(\bar v^2+ \bar v\, \bar u+5\, \bar u^2\right) \mid y_{2^{2\,\ell+1}}\end{equation}
for some non-trivial solution $\bar r, \bar s, \bar v, \bar u$ of \eqref{eq:rulMall19}, thus such that
\begin{equation}\label{NonTrivialSolution}
171\,\left( \bar r^2+ \bar  r\: \bar s\,+5\, \bar s^2\right)^2-169\,\left( \bar  v^2+ \bar  v\, \bar u+5\, \bar  u^2\right)^2 = 2 \,.
\end{equation}

Now consider the system
\begin{equation}\label{system}
\left\{
\begin{array}{lcrcl}
X &+& 171\, Y &=& \bar v^2+ \bar v\, \bar u+ 5\: \bar u^2\,, \\
X &+& 169\, Y &=& \bar r^2+ \bar r\: \bar s\,+5\, \bar s^2
\end{array}
\right.
\end{equation}
whose solution is
\[
\left\{
\begin{array}{rcl}
\overbar{X} &=& \frac{1}{2}\,\left[171\,\left(\bar r^2+ \bar r\: \bar s+5\, \bar s^2\right) - 169\,
\left(\bar v^2+ \bar v\, \bar u+ 5\, \bar u^2\right) \right]\,,
\\[.12cm]
\overbar{Y} &=& \frac{1}{2}\,\left[\left(\bar v^2+ \bar v\, \bar u+ 5\, \bar u^2\right) - \,\left(\bar r^2+ \bar r\: \bar s\,+5\, \bar s^2\right)\right]\,.\end{array}
\right.
\]

From (\ref{NonTrivialSolution}), $\bar r^2+ \bar  r\: \bar s\,+5\, \bar s^2$ and $\bar  v^2+ \bar  v\, \bar u+5\, \bar  u^2$ have the same parity; therefore, the easily checked fact that each non-trivial
integer solution $\bar r,\bar s,\bar v, \bar u$ to \eqref{NonTrivialSolution}---alias \eqref{eq:rulMall19}---satisfies $169\,(\bar r^2+\bar r\,\bar s+5\,\bar s^2)< 169\,(\bar v^2+\bar v\,\bar u+5\,\bar u^2)<171\,(\bar r^2+\bar r\,\bar s+5\,\bar s^2)$ entails that $\bar X$ and $\bar Y$ are positive integers.

 From (\ref{system}) and (\ref{NonTrivialSolution}) we get $171\,\left(\bar{X} + 169\,\bar{Y}\right)^{2} - 169\,\left(\bar{X} + 171\,\bar{Y}\right)^{2} = 2$, which simplifies into $\bar{X}^{2} - 171\cdot169\, \bar{Y}^{2} =1$. Since $\bar{Y} \neq 0$, the latter equation yields $\bar{X} = x_{g}$ and $\bar{Y} = \frac{y_{g}}{39}$, for some $g \geq  1$.
Therefore, from (\ref{system}) and (\ref{contrHyp}) we get $39\,\big(x_{g} + 169\,\frac{y_{g}}{39}\big)\big(x_{g} + 171\,\frac{y_{g}}{39}\big) \mid y_{2^{2\,\ell+1}}$, 
i.e. $\big(3\,x_{g} + 13\,y_{g}\big)\big(13\,x_{g} + 57\,y_{g}\big) \mid y_{2^{2\,\ell+1}}$,
i.e. $y_{2\,g+1} \mid y_{2^{2\,\ell+1}}$; this in turn yields $2\,g+1 \mid 2^{2\,\ell+1}$, a contradiction.
\end{proof}

\begin{corollary} $\mathcal{H}$ implies that 
the relation $\mathcal{J}_{19}(p,q)$, as defined in the Introduction, admits a polynomial Diophantine representation.
\end{corollary}
\proof[Sketch]{In close analogy with the treatment of the Heegner number $3$ as provided in \cite{CO21},
it can be shown that 
\[\begin{array}{rcl}\mathcal{J}_{19}(u,v)&\mybimp&\big[\,v\in \hug{y_{2^{2\,\ell+1}} \wt \ell \geq 0 }\setminus\hug{y_{2^{2\,\ell+1}}\wt \ell \leq 5}\,\big]\ \mywedge\ \exists\,x\big[(2\,x+1)\cdot u=v\big]\end{array}\]

\noindent is a polynomial Diophantine representation of $\mathcal{J}_{19}$\,. A salient remark is that when $b\neq 0$ the condition $\exists\,x\big[(2\,x+1)\cdot a=b\big]$ means ``$a\mid b$ and $a$ is divisible by any power of $2$ that divides $b$''; thus, since $2^{2\,\ell+2}$ is the largest power of $2$ dividing $y_{2^{2\,\ell+1}}$, in the case at hand
$\exists\,x\big[(2\,x+1)\cdot u=v\big]$ amounts to $u\mid v\ \mywedge\ u \geq 2^{2\,\ell+2}$\,.\qed
}

\medskip

Also the proof that $\mathcal{J}_{19}$ satisfies the exponential-growth criteria which are
recalled in the Introduction parallels the treatment of the Heegner number $3$ as provided in \cite{CO21}.

\section{$\mathcal{J}_{d}$ conforms to the property \eqref{propertyByMatiyasevich}: proof for the emblematic case} \label{sec:propertyByMatiyasevich}
As announced in the Introduction, Matiyasevich's property \eqref{propertyByMatiyasevich}
is satisfied by each $\mathcal{J}_d$ with $d\in\hug{2,3,7,11,19,43}$. The proofs are akin to one another for all $d$'s, just slightly simpler when $d\in\hug{2,7}$; below we illustrate how they proceed, by again choosing $d=19$ as our emblematic case-study. We will take the following two facts for granted, one of which already entered in the proof that $\mathcal{J}_{19}$ has an exponential growth while the other is proved as \cite[Lemma~6.9]{CO21}:
\begin{fact}\label{fact:one} For $n>1$, it holds that $170^{n-1}<y_n$\:. 
\end{fact}
\begin{fact}\label{fact:two} For every real number $x\geq 1$, some positive \emph{even} integer lies in the open interval $I_x\defAs\;]\,1+\log_2(1+x)\,,\:5+2\log_2 x\,[\:$. 
\end{fact}

\begin{theorem}\label{thm:matiyasevichGrowth} Take $d=19$.
Let $\alpha=170^4\:$, $\beta=2\:$, $\gamma=2^{15}\:$, and $\delta=170^3$\,. Then, to each positive integer $w$ there correspond non-negative integers $u,v$ such that
\begin{equation*}\mathcal{J}_{d}(u,v)\ \mywedge\ u<\gamma\,w^\beta\ \mywedge\ v>\delta\,\alpha^w\:.\end{equation*}
\end{theorem}
\begin{proof}
Recalling that
\[\begin{array}{rcl}\mathcal{J}_{19}(u,v)&\defFoAs&\exists\,\ell\:\left[\,\ell>5\ \mywedge\ v=y_{2^{2\,\ell+1}}\ \mywedge\ u\mid v\ \mywedge\ u\geq 2^{2\,\ell+2}\,\right],\end{array}
\]
in order to prove the claim it suffices to show that, for each $w\geq 1$, there is an integer $\ell>5$ such that
\begin{equation*}2^{2\,\ell+2}<\gamma\,w^\beta\ \ \mywedge\ \ y_{2^{2\,\ell+1}}>\delta\,\alpha^w\:.\end{equation*}

Note that \Fact{fact:one} yields $y_{2^{2\,\ell+1}}\geq 170^{2^{2\,\ell+1}-1}$ for every $\ell\in\nats$. Our task hence further reduces to proving that, for each $w\geq 1$, there is an $\ell>5$ such that
\begin{equation*}2^{2\,\ell+2}<\gamma\,w^\beta\ \ \mywedge\ \ 170^{2^{2\,\ell+1}-1}>\delta\,\alpha^w\:,\end{equation*}
namely
\begin{equation*}2^{2\,\ell+2}<2^{15}\,w^2\ \ \mywedge\ \ 170^{2^{2\,\ell+1}-1}>170^3\cdot 170^{4\,w}\:,\end{equation*}
i.e.,
\begin{equation}\label{eq:matiyasevichGrowth}2^{2\,\ell}<2^{13}\,w^2\ \ \mywedge\ \ 170^{2^{2\,\ell+1}}>170^4\cdot 170^{4\,w}\:.\end{equation}

For $w\geq 1$, \eqref{eq:matiyasevichGrowth} amounts to
\begin{equation}\label{eq:matiyasevichGrowth2}{2\,\ell}<13+2\,\log_2\,w\ \ \mywedge\ \ {2^{2\,\ell+1}}>4+{4\,w}\:.\end{equation}
In turn, ${2^{2\,\ell+1}}>4+{4\,w}$ amounts to
\begin{equation*}{2\,\ell}>1+\log_2(1+w)\:,\end{equation*}
and hence \eqref{eq:matiyasevichGrowth2} is equivalent to
\begin{equation}\label{eq:matiyasevichGrowth3}1+\log_2(1+w)<2\,\ell<13+2\,\log_2\,w\:.\end{equation}
Summing up, to get the desired claim it suffices to show that, for every $w\geq 1$, there is an $\ell>5$ satisfying
\eqref{eq:matiyasevichGrowth3}. By \Fact{fact:two}, some even integer $2\,\ell'$ satisfies $1+\log_2(1+w)<2\,\ell'<5+2\,\log_2\,w$\,. If $\ell'>5$, take $\ell=\ell'$. Otherwise take $\ell=6$: then $2\,\ell=12<13+2\,\log_2 w$, and $1+\log_2(1+w)<2\,\ell'\leq10<12$.
\end{proof}

\begin{remark} We have adopted $d=19$ as our main case-study throughout
this section, because this case is comparatively challenging; however, it is
the Heegner number $d=2$ the one leading to the candidate rule-them-all
equation that looked most promising to us when this note was first written (but see \Sec{sec:postscript}). This is why we deem it significant
to indicate which changes to the proof of \Thm{thm:matiyasevichGrowth} lead to
a proof of the analogous statement.
\begin{theorem} Take $d=2$.
Let $\alpha=3\:$, $\beta=2\:$, $\gamma=2^7\:$, and $\delta=1$\,. Then, to each positive integer $w$ there correspond non-negative integers $u,v$ such that
\begin{equation*}\mathcal{J}_{d}(u,v)\ \mywedge\ u<\gamma\,w^\beta\ \mywedge\ v>\delta\,\alpha^w\:.\end{equation*}
\end{theorem}
\proof[Sketch]{Now we must refer to the definition:
\[\begin{array}{rcl}\mathcal{J}_{2}(u,v)&\defFoAs&\exists\,\ell\,\Big[\,\ell>4\ \mywedge\ 
v={y}_{2^{\ell}}\ \mywedge\ u\mid v\ \mywedge\ u\geq 2^{\ell+1}
\,\Big].\end{array}
\]
Our previous formulations of \Fact{fact:one} can be superseded by the following:
\begin{itemize}
\item For every positive integer $n$, it holds that $3^{n-1}<y_n$ (whence
 $y_{2^\ell}\geq 3^{2^\ell-1}$).
\end{itemize}

Our goal becomes the one of proving, for each $w\geq 1$, the existence of an $\ell>4$ such that
either the condition \begin{equation*}2^{\ell+1}<\gamma\,w^\beta\ \ \mywedge\ \ y_{2^\ell}>\delta\,\alpha^w\:,\end{equation*}
or the stronger condition
\begin{equation*}2^{\ell+1}<\gamma\,w^\beta\ \ \mywedge\ \ 3^{2^\ell-1}>\delta\,\alpha^w\:,\end{equation*}
holds. The latter condition, namely $2^{\ell+1}<2^7\,w^2\ \ \mywedge\ \ 3^{2^\ell-1}>3^w$\,, gets rewritten as $2^{\ell}<2^6\,w^2\ \ \mywedge\ \ 3^{2^\ell}>3^{w+1}$\,, then as
$\ell<6+2\,\log_2\,w\ \ \mywedge\ \ {2^\ell}>1+w$\,, and then as
$\log_2(1+w)<\ell<6+2\,\log_2\,w$\,. \Fact{fact:two} yields an $\ell'$ with $\log_2(1+w)<1+\log_2(1+w)<\ell'<5+2\,\log_2\,w$; then $\ell=\max(\ell',5)$ satisfies the requirement $\ell>4$ as well, since $\ell'\leq4$ entails $\log_2(1+w)<3<5<6+2\,\log_2 w$\,.} 
\qed\end{remark}

\section{A vain attempt to find a non-trivial solution to $2\cdot\left({r}^{2} + 2\: s^{2}\right)^{2} - \left({u}^{2} + 2\: v^{2}\right)^{2}=1$}\label{sec:attemptToFindNonTrivialFor2}
Let us associate the Pell-like equation \begin{equation}\label{eq:quadratic2}2\,A^2-B^2=1\end{equation} in the unknowns $A$ and $B$ with the quaternary quartic equation corresponding to the Heegner number $d=2$\,. The solutions to this Pell-like equation in \nats form the sequence $\tup{A_0,B_0},\tup{A_1,B_1},\tup{A_2,B_2},\dots$\:, where
\begin{equation}\label{eq:quadratic2recurrence}
\begin{array}{rcrclccrclcrcl}A_0&=&B_0&=&1&\mbox{and}&\phantom{xxx}
A_{n+1}&=&3\,A_n+2\,B_n\,,&\phantom{xxx}&B_{n+1}&=&4\,A_n+3\,B_n\,,\end{array}\end{equation}
for all $n\in\nats$\,.
Our goal here is to find an integer $n>0$ such that both of
$A_{n}$ and $B_{n}$ are representable in the form $w^2+2\,t^2$\,.
We will try to achieve this by means of the method devised by Shanks and Wagstaff \cite{SW95}.

  The recurrence formula \eqref{eq:quadratic2recurrence} gives us: $A_1=5\,,\:B_1=7\,,\:A_2=29\,,\:B_2=41\,,\:A_3=169\,,\:B_3=239$\,, etc. We are less interested in the explicit values of $A_n\,,\:B_n$ than in their residue classes modulo 8; hence we form the table
  \[\begin{array}{r|ccccccccl}
    n \:&\, \text{\small1} \:&\, \text{\small2} \:&\, \text{\small3} \:&\, \text{\small4} \:&\, \text{\small5} \:&\, \text{\small6} \:&\, \text{\small7} \:&\, \text{\small8} \:&\, \dots\\\midrule
  A_n \:&\, 5 \:&\, 5 \:&\, 1 \:&\, 1 \:&\, 5 \:&\, 5 \:&\, 1 \:&\, 1 \:&\, \dots\\\midrule
  B_n \:&\, 7 \:&\, 1 \:&\, 7 \:&\, 1 \:&\, 7 \:&\, 1 \:&\, 7 \:&\, 1 \:&\, \dots\,,\end{array}\]
namely:
\[\begin{array}{lcccl}
A_n\equiv5\mod 8&&&\mbox{if }&n\equiv1,2 \mod4\:,\\
A_n\equiv1\mod 8&&&\mbox{if }&n\equiv3,0 \mod4\:,\\ 
B_n\equiv7\mod 8&&&\mbox{if }&n\mbox{ is odd}\,,\\
B_n\equiv1\mod 8&&&\mbox{if }&n\mbox{ is even}\,. 
\end{array}\]

In view of the information about representability in $\ints\left[\,\sqrt{-2}\:\right]$, as provided in \Sec{sec:uniqueFactorizationRings} (see, in particular, \Exa{exa:reprMod2}),
we only need to survey those pairs $A_n,B_n$ whose subscript $n$ is
a multiple of 4. It turns out that for no
$n\in\hug{4,8,12,16,\dots,100}$ both of $A_n$ and $B_n$ are representable. In particular:\footnote{By $p^k\|m$\,, where $p$ is a prime number and $k,m\!\in\!\nats$\,, one states that $p^k\mid m\ \mywedge\ p^{k+1}\nmid m$\,.}
\begin{description}
\item $5\,\|\,A_4=985=5\cdot 197$\,, where $5$ is inert, and so $A_4$ is not representable;
\item $103\,\|\,B_8=103\cdot 15607$\,, where $103$ is inert, and so $B_8$ is not representable;
\item $29\,\|\,A_{12}=29\cdot1549\cdot29201$\,, where $29\equiv 5\mod{8}$, hence $A_{12}$ is not representable;
\item $5\,\|\,A_{16}=5\cdot 5741\cdot52734529$\,, and hence $A_{16}$ is not representable;
\item $302663\,\|\,B_{20}=2297\cdot302663\cdot3553471$\,, where 
$302663\equiv 7\mod{8}$, hence $B_{20}$ is not representable;
\item and so on.
\end{description}

This empirical approach to finding a pair $\tup{A_n,B_n}$ of representable numbers that solves \eqref{eq:quadratic2} has hence failed; but it cannot be excluded that a larger $n$---even, maybe, infinitely many $n$'s---would work.
 
\section{P.S.: Successful calculations of non-trivial solutions to $2\cdot\left({r}^{2} + 2\: s^{2}\right)^{2} - \left({u}^{2} + 2\: v^{2}\right)^{2}=1$}\label{sec:postscript}

Readily after this note appeared on arXiv,  Evan O'Dorney (University of Notre Dame) and Bogdan Grechuk (University of Leicester) sent us kind communications that they had found non-trivial solutions to the most tantalizing quaternary quartic at stake, namely the equation $2\cdot({r}^{2} + 2\: s^{2})^{2} - ({u}^{2} + 2\: v^{2})^{2}=1$ discussed just above. On March 9, Evan O'Dorney sent us a Sage program by means of which he had located two non-trivial solutions corresponding to the Pell pairs $\tup{A_{128},B_{128}}$ and $\tup{A_{140},B_{140}}$\,. Independently of him, the next day Bogdan Grechuk sent us the following explicit values of three solutions:

\[\begin{array}{rcl}
r_1 &=& 8778587058534206806292620008143660818426865514367,\\
s_1 &=& 1797139324882565197548134105090153037130149943440,\\
u_1 &=& 5221618295817678692343699483662704959631052331713,\\ 
v_1 &=& 6739958317343073985310999451965479560858521871624;\\~\\
r_2 &=& 236514273578291664435175687910940947997062625569350147,\\ 
s_2 &=& 190287799713845710242676318005133890540427682774721600,\\ 
u_2 &=& 320623735768998122027997700001721820015837211029746377,\\ 
v_2 &=& 198400977912717981475493948031175304069461203340217424;\\~\\
r_3&=& 12046778408197306536220633533656839130825460533597900252752984198664442382891199603\_\\&&88878705520138331242183446788839967467871070821373969808509021113817760729395211791\_\\&&37051555383659,\\
s_3&=& 289694399445761398483025826848947500268215847188070555715114661420251078144498958580\_\\&&201741377919581076929573276008122608702416155843247853821556904256142740849897516701\_\\&&16452606560,\\
u_3&=& 135196357398730464210306626543956422178325136022293473333456208137694554971291574397\_\\&&360535539290774456213138712173848352564854878347629973631394683117973225078027645693\_\\&&27905413177,\\
v_3&=& 106570802532265028597850159939791513955713767327960632149341703749584264199100661176\_\\&&169441237306731010447512066786172237816690147584537869559546356198844075069003058520\_\\&&448399674296.
\end{array}\]
\noindent These are associated with the said two Pell pairs and to the companion Pell numbers ${A_{468},B_{468}}$\,.

\section*{Conclusion}

The surmise that each r.e.\ set admits a single-fold polynomial Diophantine representation
was made by Yuri Matiyasevich in 1974, inside the paper where he first proved the single-fold \emph{exponential} Diophantine representability of any r.e.\ set (cf. \cite[\pag{301}]{Mat77a}). It looks promising---but is also tantalizing---that the truth of the weaker conjecture that each r.e.\ set admits a finite-fold polynomial Diophantine representation could be established by just proving that a single fourth degree Diophantine equation 
has only a finite number of integer solutions. No algorithm can tell us, for any given
Diophantine equation, whether the set of its integer solutions is finite or infinite (cf. \cite{Davis72});\footnote{In fact, given an arbitrary instance $P=0$ of Hilbert's 10$^{\rm th}$ problem and a variable $x_0$ not occurring in $P$\,, the equation $(x_0^2+1)\cdot P=0$ is equivalent to $P=0$, with $x_0$ left free since $x_0^2+1$ never vanishes; hence it has infinitely many integer solutions if and only if the equation $P=0$ has at least one integer solution. Being able to determine, for any given Diophantine equation $Q=0$, whether or not it is finite-fold, would hence enable one to solve Hilbert's 10$^{\rm th}$ problem. This would contradict the today known undecidability of the latter problem (cf. \cite{Davis73,DMR76}).} nevertheless, the structure of our six candidate rule-them-all equations is so simple that we may hope that the finite-foldness of at least one of them will come to light.

It must be mentioned that, in \cite{Tyszka18}, Apoloniusz Tyszka advances a conjecture that,
if true, would falsify Matiyasevich's finite-fold representability conjecture, which has been,
throughout, our polar star.

\bigskip

\noindent{\bf Acknowledgements.}\ %
We gratefully acknowledge partial support from project ``STORAGE---Uni\-versit\`{a} degli Studi di Catania, Piano della Ricerca 2020/2022, Linea di intervento 2'', and from INdAM-GNCS 2019 and 2020 research funds.

\medskip

Discussions with Pietro Corvaja were very fruitful for the matters of this paper.

%
%
%
\bibliographystyle{alpha}

\end{document}